\documentclass{amsart}

\usepackage{amsmath,amssymb,amsthm}
\usepackage{graphicx,tikz}
\newtheorem{theorem}{Theorem}

\theoremstyle{definition}

\theoremstyle{remark}

\begin{document}

\title[]{A Pointwise Inequality for Derivatives\\ of Solutions of the Heat Equation\\ in Bounded Domains}
\subjclass[2010]{35J05, 35K05, 35P05} 
\keywords{Heat Equation, Laplacian eigenfunction, Hessian estimate.}
\thanks{S.S. is supported by the NSF (DMS-1763179) and the Alfred P. Sloan Foundation.}

\author[]{Stefan Steinerberger}
\address{Department of Mathematics, University of Washington, Seattle}
\email{steinerb@uw.edu}

\begin{abstract} Let $u(t,x)$ be a solution of the heat equation in $\mathbb{R}^n$. Then, each $k-$th derivative also solves the heat equation and
satisfies a maximum principle, the largest $k-$th derivative of $u(t,x)$ cannot be larger than the largest $k-$th derivative of $u(0,x)$. We
prove an analogous statement for the solution of the heat equation on bounded domains $\Omega \subset \mathbb{R}^n$ with Dirichlet boundary conditions.
As an application, we give a new and fairly elementary proof of the sharp growth of the second derivatives of Laplacian eigenfunction $-\Delta \phi_k = \lambda_k \phi_k$ with Dirichlet conditions on smooth domains $\Omega \subset \mathbb{R}^n$.
 \end{abstract}

\maketitle

\section{Introduction and Results}

In Euclidean space $\mathbb{R}^n$, it is always possible to exchange the order of partial derivatives: in particular, derivatives of solutions of the heat equation also satisfy the heat equation and therefore enjoy many nice properties, in particular a maximum principle. 
Let $u(t,x)$ denote a solution of the heat equation in $\mathbb{R}^n$. Then it is explicitly given by
$$ u(t,x) = \frac{1}{(4\pi t)^{d/2}} \int_{\mathbb{R}^n} \exp\left(-\frac{\|x-y\|^2}{4t} \right) u(0,y) dy.$$
At any given point $x_0 \in \mathbb{R}^2$ and any unit vector $\nu \in \mathbb{S}^{n-1}$ and any $k \geq 1$, we can (at least formally) differentiate under the integral sign to obtain
$$ \frac{\partial^k u}{\partial \nu^k} (t, x_0) = \frac{1}{(4\pi t)^{d/2}} \int_{\mathbb{R}^n}\exp\left(-\frac{\|x_0-y\|^2}{4t} \right)\frac{\partial^k u}{\partial \nu^k} (0, x_0)dy.$$
 We were interested in whether there is an analogous result on bounded domains in terms of the heat kernel $p_t(\cdot, \cdot)$ of the domain $\Omega$.

\begin{theorem}[Main Result] Let $\Omega \subset \mathbb{R}^n$ be a domain with $C^1-$boundary and let $f \in C^k(\Omega) \cap L^{\infty}(\Omega)$, let $x_0 \in \Omega$ and $\nu \in \mathbb{S}^{n-1}$. Then, for all $t>0$, all , the solution $e^{t\Delta} f$ of the heat equation with Dirichlet boundary conditions $u(t, x) = g(x)$ has 
$$X =  \left| \frac{\partial^k }{\partial \nu^k} e^{t\Delta} f(x_0) - \int_{\Omega}   p_t(x_0, y) \frac{\partial^k }{\partial \nu^k} f(y) dy \right|$$
bounded by
$$X \leq  \left(1 - \int_{\Omega}p_t(x_0, y) dy \right) \max_{0 \leq s \leq t}  \left\| \frac{\partial^k }{\partial \nu^k} e^{s\Delta} f \right\|_{L^{\infty}(\partial \Omega)}.$$\\
\end{theorem}

 Note that if $t \ll d(x_0, \partial \Omega)^2$, then the integral term in our error bound is actually quite small independently of what happens on the boundary. 
In particular, in free space, the heat kernel always has total integral 1 and we get $X \equiv 0$. 
  The result seems to have fairly natural extensions to the Neumann Laplacian, higher derivatives and even more general parabolic equations following essentially the same type of argument.
  This type of argument might have an interesting analogue on manifolds (both with and without boundary). In the case of manifolds without boundary, one would expect that the underlying curvature has an additional perturbative effect on the particles (`stochastic parallel transport', see e.g. Bismut \cite{bis}, Elworthy \& Li \cite{el} and  Thalmaier \& Wang \cite{thal}). We also refer to recent developments on second-order Feynman Kac formulas (see Li \cite{li} and Thompson \cite{thomp}).\\

\textbf{An Application.}
 Let $\Omega \subset \mathbb{R}^n$ be a bounded domain with smooth boundary. We consider Laplacian eigenfunctions,
$ -\Delta \phi_k = \lambda_k \phi_k$,
with Dirichlet conditions on the boundary $\partial \Omega$. If these eigenfunctions to be normalized in $L^2$, i.e. $\| \phi_k\|_{L^2}=1$, then results of
Levitan in 1952 \cite{levitan}, Avakumovic in 1956 \cite{ava} and H\"ormander in 1968 \cite{hor} guarantee that
$$ \left\| \phi_k\right\|_{L^{\infty}(\Omega)} \leq c_{\Omega} \cdot \lambda_k^{\frac{n-1}{4}},$$
where $c_{\Omega}$ is a constant depending only on $\Omega$.
This estimate is sharp (for an example on a ball, see \cite[\S 2.3]{grieser}) and has been well studied \cite{berard, blair, hassell, smith0, smith, sogge, sogge2, sogge3, sogge35, sogge4}.

The optimal estimate for the gradient is
$$ \left\| \nabla \phi_k\right\|_{L^{\infty}(\Omega)} \leq c_{\Omega} \cdot \lambda_k^{\frac{n+1}{4}}$$
and has been studied by Xu \cite{xu1, xu2, xu3} and, subsequently, by
 Arnaudon, Thalmaier \& Wang \cite{arnaud}, Cheng, Thalmaier \& Thompson \cite{cheng}, Hu, Shi \& Xu \cite{hu} and Shi \& Xu \cite{shi}. On compact manifolds without boundary, there are results of Xu \cite{xu0} for all derivatives and by Wang \& Zhou \cite{wang} for linear combinations of eigenfunctions. Recently, Frank \& Seiringer \cite{frank} showed that for compact $\Omega \subset \mathbb{R}^n$ with $C^{k,\delta}-$smooth boundary, there is an estimate
 $$ \left\| \frac{\partial^k \phi_k}{\partial \nu^k} \right\|_{L^{\infty}(\Omega)} \lesssim_{\Omega, k}~  \lambda^{k/2} \| \phi_k\|_{L^{\infty}(\Omega)} $$
As an application, we give an elementary proof of this inequality for $k=2$.
\begin{theorem} Let $\Omega \subset \mathbb{R}^n$ be a compact domain with smooth boundary. There exists a constant $c_{\Omega}$ such that for all solutions of $-\Delta \phi_k = \lambda_k \phi_k$ that vanish on $\partial \Omega$
$$ \| D^2 \phi_k(x)\|_{L^{\infty}(\Omega)} \leq c_{\Omega} \cdot  \lambda_k^{1/2} \cdot \| \nabla \phi_k\|_{L^{\infty}(\Omega)}.$$
\end{theorem}
More generally, Theorem 1 allows us to obtain similar estimates also for higher derivatives provided there is some control on the derivatives on the boundary.

\section{Proof of Theorem 1}

\begin{proof} We describe the proof for $k=2$ in detail: the more general case $k \geq 3$ is completely analogous (after replacing the second differential quotient by the corresponding $k-$th differential quotient). We use a probabilistic argument. For any $x \in \Omega$, let $\omega_x(t)$ denote a Brownian motion started in $x$ after $t$ units of time. This Brownian motion gets `stuck' once it hits the boundary (this corresponds to Dirichlet boundary conditions). This gives us a way of solving the heat equation via
$$ e^{t\Delta} f(x) = \mathbb{E} \left(f(\omega_x(t))\right),$$
where the expectation ranges over all Brownian motions $\omega_x$ started in $x$ that run for $t$ units of time.
Our goal is to control the size of second derivatives of $e^{t\Delta} f(x)$ for points inside the domain. Let now $x_0 \in \Omega$ be fixed and let $\nu \in \mathbb{S}^{n-1}$ be some fixed direction. Calculus tells us that the second derivative of $e^{t\Delta} f$ at $x_0$ in direction $\nu$ is given by the limit of a differential quotient
$$  \left(\frac{\partial^2 }{\partial \nu^2}  e^{t\Delta} f \right)(x_0) = \lim_{\varepsilon \rightarrow 0} \frac{ e^{t\Delta} f(x_0 + \varepsilon \nu) - 2 e^{t\Delta} f(x_0) + e^{t\Delta} f(x_0 - \varepsilon \nu)}{\varepsilon^2}.$$
We will, for the remainder of the proof, control exactly this differential quotient by carefully grouping the three expectations arising from
\begin{align*}
e^{t\Delta} f(x_0 + \varepsilon \nu) - 2 e^{t\Delta} f(x_0) + e^{t\Delta} f(x_0 - \varepsilon \nu) &= \mathbb{E} \left( f(\omega_{x+ \varepsilon \nu}(t))\right) - 2 \cdot \mathbb{E} \left( f(\omega_{x}(t))\right) \\
&+\mathbb{E} \left(f(\omega_{x- \varepsilon \nu}(t))\right).
\end{align*}
These are three independent Brownian motions started in three different points. However, these three initial points are very close to one another (ultimately, $\varepsilon \rightarrow 0$), so we expect them to be somewhat related. Let $A$ denote all Brownian motion paths started in $0 \in \mathbb{R}^n$ and running for $t$ units of time. Then, for each $y \in \Omega$, we can use translation invariance of Brownian motion in $\mathbb{R}^n$ to write the expectations as an expectation over the set $A$ via
$$ e^{t\Delta} f(y) = \mathbb{E} \left(f(\omega_y(t))\right) =  \mathbb{E}_{a \in A} \left( f(a(t) + y) \cdot 1_{\left\{a(s) + y \in \Omega ~\mbox{\tiny for all}~0 \leq s \leq t\right\}}\right).$$
This has the advantage of being able to take the expectation with respect to one universal set $A$ shared by all three Brownian motions.
The size of the characteristic function has an analytic expression which is given by the heat kernel $p_t(\cdot, \cdot)$: 
$$ \mathbb{P} \left( 1_{\left\{a(s) + y \in \Omega ~\mbox{\tiny for all}~0 \leq s \leq t\right\}} \right) = \int_{\Omega} p_t(y, z) dz.$$
There is an interesting subset of $A$ depending on $x_0, \nu$ and $\varepsilon$
 $$ A_{\varepsilon} = \left\{ a \in A: \quad \begin{cases} a(s) + x_0 \in \Omega \hspace{5pt} \qquad \mbox{\tiny for all}~0 \leq s \leq t \\
  a(s) + x_0 + \varepsilon \nu \in \Omega ~\mbox{\tiny for all}~0 \leq s \leq t \\
   a(s) + x_0 - \varepsilon \nu \in \Omega ~\mbox{\tiny for all}~0 \leq s \leq t  \end{cases}\right\} .$$
$A_{\varepsilon} \subset A$ is the set of Brownian paths that remain in $\Omega$ for all time $0\leq s \leq t$ independently in which of the three points $x_0, x_0 \pm \varepsilon \nu$ they are started in.  For paths in $A_{\varepsilon}$, the differential quotient is easy to analyze: the path ends at a certain point $a(t)$ and the relative position of the three Brownian particles has been preserved since not a single one of them has hit the boundary. Recalling that for $B \subset \Omega$
$$ \mathbb{P}\left(\omega_x(t) \in B\right) = \int_{B} p_t(x_0, y) dy,$$
we see that
\begin{align*}
 X = \lim_{\varepsilon \rightarrow 0} \frac{1}{\varepsilon^2}\mathbb{E}_{A_{\varepsilon}} \left[ f(\omega_{x_0+ \varepsilon v}(t)) - 2 f(\omega_{x_0}(t)) + f(\omega_{x_0- \varepsilon v}(t)) \right]   \end{align*}
 can be written as
 \begin{align*}
 X = \lim_{\varepsilon \rightarrow 0}\mathbb{E}_{a \in A_{\varepsilon}} \left( \frac{\partial^2 f}{\partial \nu^2}(a(t)) \right) = \int_{\Omega} p_t(x_0, z)  \frac{\partial^2 f}{\partial \nu^2}(z) dz.
 \end{align*}
 
This serves as a probabilistic proof that in $\mathbb{R}^n$, solving the heat kernel and differentiation commute since in that case $A_{\varepsilon} = A$ because there is no boundary. The remainder of the argument will be concerned with $A \setminus A_{\varepsilon}$. 
The cases that are in $A \setminus A_{\varepsilon}$ can be written as a disjoint union
$$ A \setminus A_{\varepsilon} = \bigcup_{0 \leq s_0 \leq t} A_{1,s_0} \cup A_{2,s_0} \cup A_{3, s_0},$$
 where $A_{i, s_0}$ is the event where the $i-$th of the three particles (having enumerated them in an arbitrary fashion) has hit the boundary at time $s_0$ and is the first of the three particles to do so (in case two particles hit the boundary simultaneously, we may put this event into either of the two sets; if all three hit simultaneously, it can be put into any of the three sets). We have a very good understanding of the size of $A \setminus A_{\varepsilon}$ since
 $$ \lim_{\varepsilon \rightarrow 0} \mathbb{P}\left(A \setminus A_{\varepsilon}\right)  = 1 - \int_{\Omega} p_t(x_0, y) dy.$$
 It remains to understand the expected value of the differential quotient conditional on the path being in $A_{i, s_0}$ for all $1 \leq i \leq 3$ and all $0 \leq s_0 \leq t$. As it turns out, these cases can all be analyzed in the same fashion.
 We illustrate the argument using the example shown in Fig. \ref{fig:conf}. In that example, the middle particle has impacted the boundary at time $s_0$. That middle Brownian motion is stuck and 
 $$ \mathbb{E} \left( f(\omega_{x_0}(t)) \big| \omega_{x_0}(s_0) = a(s_0) + x_0  \right) = 0.$$

\begin{center}
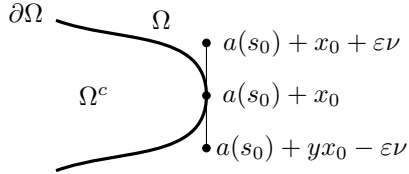
\begin{figure}[h!]
\begin{tikzpicture}
\draw [very thick] (3,0) to[out=20, in = 270]  (5,1) to[out=90, in =340] (3,2);
\filldraw (5,0.3) circle (0.05cm); 
\node at (6, 1) {$a(s_0) + x_0$};
\node at (4.4, 2) {$\Omega$};
\node at (3.5, 1) {$\Omega^c$};
\node at (2.6, 2.1) {$\partial \Omega$};
\node at (6.4, 1.7) {$a(s_0) + x_0 + \varepsilon \nu$};
\node at (6.4, 0.3) {$a(s_0) + yx_0- \varepsilon \nu$};
\filldraw (5,1) circle (0.05cm);
\filldraw (5,1.7) circle (0.05cm);
\draw (5, 0.3) -- (5, 1.7);
\end{tikzpicture}
\caption{A path in $ A_{2, s_0} \subset A \setminus A_{\varepsilon}$: the middle point hits the boundary.}
\label{fig:conf}
\end{figure}
\end{center}

It remains to understand the expected value of $f(\omega_{x_0}(t))$ subject to knowing that at time $s_0$ the particle is in the position $a(s_0) + x_0 + \varepsilon \nu$ or $a(s_0) + x_0 - \varepsilon \nu$. At this point we use Markovianity: the Brownian motion does not remember its past and behaves as if it were freshly started in that point. In particular, this shows that
$$ \mathbb{E} \left( f(\omega_{x_0}(t)) \big| \omega_{x_0}(s_0) = a(s_0) + x_0 + \varepsilon \nu \right) = \mathbb{E} \left( f(\omega_{a(s_0) + x_0 + \varepsilon \nu}(t-s_0)) \right).$$
However, this is merely the formula for the solution of the heat equation and
$$  \mathbb{E} \left( f(\omega_{a(s_0) + x_0 + \varepsilon \nu}(t-s_0)) \right) = e^{(t-s_0)\Delta} f(a(s_0) + x_0 + \varepsilon \nu).$$
Likewise, we have that
$$ \mathbb{E} \left( f(\omega_{x_0}(t)) \big| \omega_{x_0}(s_0) = a(s_0) + x_0 - \varepsilon \nu \right) = e^{(t-s_0)\Delta} f(a(s_0) + x_0 - \varepsilon \nu).$$
Finally, we argue that even the middle point (the one that already impacted on the boundary) can be written the same way since
 $$ \mathbb{E} \left( f(\omega_{x_0}(t)) \big| \omega_{x_0}(s_0) = a(s_0) + x_0  \right) = 0 =  e^{(t-s_0)\Delta} f(a(s_0) + x_0).$$

However, these three identities tell us a nice story: it tells us that evaluating the `probabilistic' second differential quotient amounts to, in this special case,
merely to evaluating the second differential quotient of
$$ e^{(t-s_0)\Delta} f \qquad \mbox{at the point}~a(s_0) + x_0\qquad \mbox{in direction}~\nu.$$
This leads to the following conclusion: for any smooth, compact domain, the heat equation started with $u(0,x) = f(x)$ inside $\Omega$ and constant boundary conditions $u(t,x) = g(x)$ for $x \in \partial \Omega$ has a solution
$$ u(t,x) = \int_{\Omega} p_t(x,y) f(y) dy + \int_{\partial \Omega} q_t(x,y)g(y) dy,$$
where $q_t(x,y)$ converges to the harmonic measure as $t \rightarrow \infty$ (since $\Omega$ is smooth, regularity of the boundary does not play a role).
From this we can deduce that
$$  \left(\frac{\partial^2 e^{t\Delta} f }{\partial \nu^2}  \right)(x) = \int_{\Omega} p_t(x,y) \frac{\partial^2 f}{\partial \nu^2}(y) dy + \int_{0}^{t} \int_{\partial \Omega} \frac{\partial q_s}{\partial s}(x,y)  \frac{\partial^2 e^{(t-s) \Delta} f}{\partial \nu^2}(y) dy ds.$$
 The second integral can be easily bounded from above by
 \begin{align*}
  \left|  \int_{0}^{t} \int_{\partial \Omega} \frac{\partial q_s}{\partial s}(x,y)  \frac{\partial^2 e^{(t-s) \Delta} f}{\partial \nu^2}(y) dy ds \right| &\leq
   \int_{0}^{t} \int_{\partial \Omega} \frac{\partial q_s}{\partial s}(x,y) \left| \frac{\partial^2 e^{(t-s)\Delta}f}{\partial \nu^2}(y) \right| dy ds\\
&=  \int_{\partial \Omega} q_t(x,y) \max_{0 \leq s \leq t} \left| \frac{\partial^2 e^{(t-s)\Delta} f}{\partial \nu^2}(y)  \right| dy.
 \end{align*}
 In particular, recalling that $u(t,x) \equiv 1$ is a solution of the heat equation with initial data $u(0,x) \equiv 1$ and boundary data $g(x) \equiv 1$, we have that
 $$ \int_{\Omega} p_t(x,y) dy + \int_{\partial \Omega} q_t(x,y) dy = 1.$$
 This shows that the integral
 $$ I = \int_{\partial \Omega} q_t(x,y) \max_{0 \leq s \leq t} \left| \frac{\partial^2 e^{(t-s)\Delta} f}{\partial \nu^2}(y)  \right| dy$$
 can be bounded by
 $$ I\leq  \left(1 - \int_{\Omega}p_t(x_0, y) dy \right) \max_{0 \leq s \leq t}  \max_{y \in \partial \Omega} \left|  \frac{\partial^2 e^{(t-s)\Delta} f}{\partial \nu^2}  \right|,$$
 where $\nu = \pm \nu$ always points inside the domain. This is the desired statement. \end{proof}
We note that the last few steps are certainly a bit wasteful and one could obtain more precise estimates if one were to assume additional knowledge about $p_t(\cdot, \cdot)$, $q_t(\cdot, \cdot)$ or the harmonic measure $\omega_{x_0}$.

\section{Proof of Theorem 2}
The proof decouples into the following steps. 
\begin{enumerate}
\item First, we show that second derivatives on the boundary are controlled and at most of size $\lesssim_{\partial \Omega} \| \nabla \phi_k\|_{L^{\infty}}$. The implicit constant will only depend on the mean curvature of the boundary $\partial \Omega$.
\item We then use this in combination with Theorem 1. We use time scale $t = \varepsilon \lambda_k^{-1}$ and will show that the entire argument
can be carried out with a sufficiently small $\varepsilon>0$ whose final size only depends on the geometry of $\partial \Omega$ and the dimension $n$. If the second derivatives assume their maximum value in a point $x_0$, then there exists a set $A$ in a $\sqrt{t}-$neighborhood of $x_0$ where the second derivatives are large. $A$ itself is large in the sense of
$$  \int_{A}   p_t(x_0, y) dy \geq 1 - 4 \varepsilon.$$
\item This implies that $x_0$ is not too close to the boundary: $d(x_0, \partial \Omega) \geq \sqrt{t}$.
\item Finally, we show that if we consider a ball of radius $\sqrt{t}$ around $x_0$ (and, by the previous step, this ball is fully contained in $\Omega$), then there exists a line segment such that second derivatives are large on most of its length.
\item The Fundamental Theorem of Calculus then shows that the derivatives have to grow very quickly along the line segment and this will lead to a contradiction once the second derivatives are too large: this will show that $\| \nabla \phi_k\|_{L^{\infty}}$ has to be large which then contradicts known bounds.
\end{enumerate}

\subsection{Eigenfunctions on the boundary.} 
We start the argument by noting that eigenfunctions $-\Delta \phi_k = \lambda_k \phi_k$ in smooth domains with vanishing Dirichlet boundary conditions on $\partial \Omega$ cannot have a particularly large second derivative on the boundary. We will prove this by expressing the Laplacian in local coordinates at the boundary: more precisely, let $\partial \Omega$ be a hypersurface in $\mathbb{R}^n$ and let $\nu$ be a unit normal vector to $\partial \Omega$. Using $\Delta_{\partial \Omega}$ to denote the Laplacian in the induced metric on $\partial \Omega$ and $\Delta$ to denote the classical Laplacian in $\mathbb{R}^n$, we have the identity
$$ \Delta u = \Delta_{\partial \Omega} u + (n-1)H \frac{\partial u}{\partial n} + \frac{\partial^2 u}{\partial n^2},$$
where $H$ is the mean curvature of the boundary $\partial \Omega$ in that point. This identity is derived, for example, in the book by Sperb \cite[Eq. 4.68]{sperb}. 

\begin{center}
\begin{figure}[h!]
\begin{tikzpicture}
\draw [thick] (0,0) to[out= 320, in =190] (4, -1);
\filldraw (2,-1.02) circle (0.04cm);
\draw [->] (2, -1.02) -- (2.18, -0.2);
\node at (1.9, -1.3) {$x$};
\node at (2.4, -0.3) {$\nu$};
\node at (4.5, -1) {$\partial \Omega$};
\end{tikzpicture}
\caption{A point on the boundary and a normal direction.}
\end{figure}
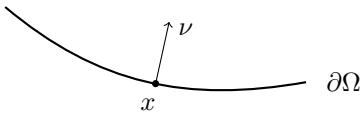
\end{center}

This equation is particularly useful in our case because $\phi_k$ vanishes identically on the boundary and thus $\Delta_{\partial \Omega} \phi_k = 0$. The identity is also frequently used on level sets of $u$ (since then the Laplacian $\Delta_{\partial \Omega} u$ vanishes for the same reason), see for example Kawohl \& Horak \cite{kawohl}. We refer to Reilly \cite{reilly} for a friendly introduction to this identity in low dimensions. Since $\phi_k$ vanishes on the boundary, we have
$$  \frac{\partial^2 \phi_k}{\partial \nu^2} + (n-1) H \frac{\partial \phi_k}{\partial \nu} = 0.$$

However, the mean curvature $H$ is bounded depending only on $\partial \Omega$ because the domain is smooth. Therefore, using the gradient estimate of Hu, Shi \& Xu \cite{hu} on compact Riemannian manifolds with boundary, we get
\begin{align*}
 \left| \frac{\partial^2 \phi_k}{\partial \nu^2}  \right| &\leq n \| H\|_{L^{\infty}(\partial \Omega)} \left| \frac{\partial \phi_k}{\partial \nu}  \right| \lesssim_{\Omega} \| \nabla \phi_k\|_{L^{\infty}} \lesssim_{\Omega} \lambda^{1/2} \| \phi_k\|_{L^{\infty}} \lesssim \lambda_k^{\frac{n+1}{4}}.
 \end{align*}
This shows that second derivatives in the normal direction cannot be much larger than first derivatives, since they are generated as a combination of first derivatives and the curvature. As for second derivatives in directions orthogonal to the normal direction, we note that due to the smoothness of $\Omega$, points at distance $\varepsilon$ are $\sim \varepsilon^2$ away from the boundary, where the implicit constant depends on the curvature of $\partial \Omega$. This allows us to bound the second differential quotient by  $\leq c_{\Omega} \| \nabla f\|_{L^{\infty}},$
where $c_{\Omega}$ depends only on the local curvature of $\partial \Omega$. 

\subsection{Applying Theorem 1.} Let us fix $x_0 \in \Omega$ and $\nu_0 \in \mathbb{S}^{n-1}$ so that the second derivative in $x_0$ in direction $\nu_0$ is among the largest that can occur, i.e.
$$   \frac{\partial^2 \phi_k}{\partial \nu_0^2}(x_0)  = \max_{x \in \Omega, \nu \in \mathbb{S}^{n-1}} \left| \frac{\partial^2 \phi_k}{\partial \nu^2}(x) \right| = c_1 \lambda_k^{\frac{n+3}{4}}.$$
We assume that this largest second directional derivative is positive (for ease of exposition): if it is negative, we consider without loss of generality $-\phi_k$ instead. Our goal is to now deduce a contradiction once $c_1$ is sufficiently large.
Let $\varepsilon$ be a small parameter (our goal will be to show that there is a sufficiently small but positive parameter $\varepsilon > 0$ depending only on $\Omega$ such that all the subsequent arguments work). The first step works for all $\varepsilon>0$: we use Theorem 2 for time $t = \varepsilon \lambda_k^{-1}$ in combination with
$$ e^{t\Delta} \phi_k(x) = e^{-\lambda_k t} \phi_k(x)$$
 to obtain
$$ e^{-\varepsilon} \frac{\partial^2 }{\partial \nu_0^2}  \phi_k(x_0) \leq \int_{\Omega}   p_t(x_0, y) \frac{\partial^2 }{\partial \nu_0^2} \phi_k(y) dy   + c_2 \lambda_k^{\frac{n+1}{4}},$$
where the existence of an absolute constant $c_2$ depending only on $\Omega$ follows from \S 3.1. We introduce the set where second directional derivatives in direction $\nu_0$ are `large' (positive and at least half the value of the maximum) 
$$ A =\left\{x \in \Omega:  \frac{\partial^2\phi_k }{\partial \nu_0^2} (x) \geq  \frac{1}{2}\frac{\partial^2 \phi_k}{\partial \nu_0^2}  (x_0) \right\}.$$
This allows us to bound the inequality above by replacing the partial derivatives in $y$ by the maximum partial derivative 
\begin{align*}
e^{-\varepsilon} \frac{\partial^2 }{\partial \nu_0^2}  \phi_k(x_0) &\leq \int_{\Omega}   p_t(x_0, y) \frac{\partial^2 }{\partial \nu_0^2} \phi_k(y) dy   + c_2 \lambda_k^{\frac{n+1}{4}} \\
&\leq   \frac{1}{2} \int_{\Omega \setminus A}   p_t(x_0, y) \left(\frac{\partial^2 }{\partial \nu_0^2} \phi_k(x_0)\right) dy    \\
&+ \int_{A}   p_t(x_0, y) \left(\frac{\partial^2 }{\partial \nu_0^2} \phi_k(x_0) \right) dy + c_2 \lambda_k^{\frac{n+1}{4}}.
\end{align*}
Dividing by the largest derivative, we get
$$ 1 - \varepsilon \leq e^{-\varepsilon} \leq \frac{1}{2} \int_{\Omega \setminus A}   p_t(x_0, y) dy + \int_{A}   p_t(x_0, y) dy + \frac{c_2}{c_1} \frac{1}{\sqrt{\lambda_k}}.$$
Recall that
$$  \int_{\Omega \setminus A}   p_t(x_0, y) dy + \int_{A}   p_t(x_0, y) dy = \int_{\Omega} p_t(x,y) dy \leq 1.$$
This shows that once $\lambda_k$ is sufficiently large, say, so large that
$$ \frac{c_2}{c_1} \frac{1}{\sqrt{\lambda_k}} \leq \varepsilon$$
we have
$$ 1 - 2\varepsilon \leq \frac{1}{2} \int_{\Omega \setminus A}   p_t(x_0, y) dy + \int_{A}   p_t(x_0, y) dy.$$
However, if $1 - 2\varepsilon \leq a/2 + b$ and $a+b \leq 1$, then $a \leq 4\varepsilon$ and thus
$$ \int_{A}   p_t(x_0, y) dy \geq 1 - 4 \varepsilon.$$
This means that if we start a Brownian motion in $x_0$ and let it run for $t = \varepsilon \lambda_k^{-1}$ units of time, then the likelihood of never hitting the boundary and ending up in the set $A$ is actually quite large. We note that this argument can be used for each $\varepsilon > 0$ at the cost of excluding finitely many initial eigenfunctions (whose number depends on $\varepsilon$ in the manner outlined above).

\subsection{$x_0$ is far away from the boundary.} This section uses the inequality
$$ \int_{A}   p_t(x_0, y) dy \geq 1 - 4 \varepsilon$$
 to prove that $x_0$ is at least $\geq c_3 \cdot \sqrt{t}$ away from the boundary. Indeed, we obtain a slightly stronger result and show that we could choose $c_3 = 1$ for $\varepsilon$ sufficiently small (depending only on $\Omega$). Since the boundary is smooth and compact, there exists an length scale $\delta_0$ such that at the scale of $\delta_0$ (or below) the boundary behaves roughly like a hyperplane. Suppose $d(x_0, \partial \Omega) \leq 0.01  \sqrt{t} \ll \delta_0$. Then the geometry looks a bit as in Fig. \ref{fig:wavel}.

 \begin{center}
\begin{figure}[h!]
\begin{tikzpicture}[scale=1.3]
\draw [thick] (0,0) to[out=70, in=280] (0,3);
\filldraw (1,1.5) circle (0.04cm);
\node at (1.3, 1.5) {$x_0$};
\draw[thick] (1, 1.5) circle (1.5cm);
\draw [<->] (1, 1.6) -- (1, 2.9);
\node at (1.3, 2.4) {$\sqrt{t}$};
\node at (-0.1, 2) {$\partial \Omega$};
\draw [<->] (0.3, 1.5) -- (0.9, 1.5);
\node at (0.8, 1.2) { $d(x_0, \partial \Omega)$};
\end{tikzpicture}
\caption{A maximal point being closer than $\sqrt{t}$ to the boundary.}
\label{fig:wavel}
\end{figure}
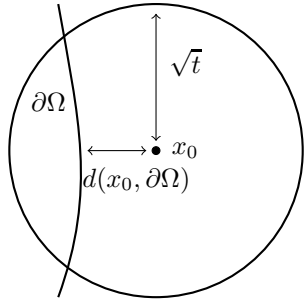
\end{center}
Since $\partial \Omega$ is smooth and we operate at small length scales (relative to $\partial \Omega$), we know that $d(x_0, \partial \Omega) \leq 0.01 \sqrt{t}$ implies that $B_{\sqrt{t}}(x_0)$ intersects the boundary $\partial \Omega$ in a large segment. This also shows that a typical Brownian motion run for $t$ units of time will hit the boundary with nontrivial probability bounded away from 0 which contradicts our lower bound on the survival probability of Brownian motion up to time $t$
$$ 1 - 4 \varepsilon \leq \int_{A}   p_t(x_0, y) dy \leq \int_{\Omega}   p_t(x_0, y) dy.$$
It remains to make this more quantitative. There is a simple way of bounding the survival probability of a Brownian motion as follows: let us assume that $x_0$ is in the origin, the boundary $\partial \Omega$ is oriented so that it is given by the hyperplane $\left\{x \in \mathbb{R}^n: x_1 = d \right\}$ at distance $d$ from the origin. The boundary $\partial \Omega$ is \textit{not} a hyperplane but for sufficiently small length scales behaves effectively as one, the curvature acts as a lower order term (which we will account for in the subsequent argument). The main question is then whether, within $t$ units of time, the first component of the Brownian motion is ever larger than $d$.
 The first component of an $n-$dimensional Brownian motion is a one-dimensional Brownian motion and we can apply the reflection principle to conclude
$$ \mathbb{P}\left( \max_{0 \leq s \leq t} B(s) \geq d \right) = \mathbb{P}\left( |B(t)| \geq d \right).$$
This refines our estimate to
\begin{align*}
 1 - 4 \varepsilon &\leq \int_{A}   p_t(x_0, y) dy \leq \int_{\Omega}   p_t(x_0, y) dy \\
 &\leq 1 - \mathbb{P}\left[ |B(t)| \geq  2 \cdot d(x_0, \partial \Omega) \right],
 \end{align*}
where the factor of 2 compensates for the higher order term coming from the curvature of $\partial \Omega$ (and thus valid for $t$ sufficiently small depending only on $\partial \Omega$). Thus
$$ \mathbb{P}\left[ |B(t)| \geq  2 \cdot d(x_0, \partial \Omega) \right] \leq 4 \varepsilon.$$
$B(t)$ is distributed like the Gaussian $\mathcal{N}(0,t)$ and thus, by rescaling,
$$ \mathbb{P}\left[ |B(t)| \geq  2 \cdot d(x_0, \partial \Omega) \right] = \mathbb{P}\left[ |B(1)| \geq  \frac{2}{\sqrt{t}} \cdot d(x_0, \partial \Omega) \right] \leq 4 \varepsilon.$$
Knowing that this quantity is less than $4\varepsilon$ will lead to a lower bound on $d(x_0, \partial \Omega)$. The standard tail bound, valid for all $z>0$,
$$ \mathbb{P} (B(0,1) \geq z) \geq \frac{1}{\sqrt{2\pi}} \left(\frac{1}{z} - \frac{1}{z^3} \right)e^{-z^2/2}$$
implies, for $z \geq 2$,
$$ 2\varepsilon  \geq \mathbb{P} (B(0,1) \geq z) \geq \frac{1}{2\pi z} e^{-z^2/2}.$$
We want to argue that this forces
$$ z \geq z_0 = \sqrt{\log{\frac{1}{\varepsilon}}}$$
because plugging in $z_0$ leads to
 $$ \frac{1}{2\pi z_0} e^{-z_0^2/2} =  \frac{1}{2\pi}  \frac{1}{\sqrt{\log{\frac{1}{\varepsilon}}}} e^{-z_0^2/2} =   \frac{1}{2\pi}  \frac{\sqrt{\varepsilon}}{\sqrt{\log{\frac{1}{\varepsilon}}}} $$
 which is larger than $2 \varepsilon$ for all $0 \leq \varepsilon \leq 0.0009$. Thus, for $\varepsilon$ sufficiently small (which here is an absolute constant), we have the desired inequality.
This forces 
$$z \geq \sqrt{\log{\frac{1}{\varepsilon}}} $$
 and thus
$$ d(x_0, \partial \Omega) \geq \frac{ \sqrt{t}}{4} \cdot \sqrt{\log{\frac{1}{\varepsilon}}}. $$
Recalling that $t = \varepsilon \lambda_k^{-1}$, we have
$$ d(x_0, \partial \Omega) \geq \frac{1}{4} \cdot \sqrt{\varepsilon} \cdot \sqrt{\log{\frac{1}{\varepsilon}}} \cdot \lambda_k^{-1/2}.$$
 For $\varepsilon \leq e^{-16}$, it is at least $\varepsilon^{1/2} \lambda_k^{-1} = \sqrt{t}$ away from the boundary. We will not, strictly speaking, need this and could absorb any constant that arises here in the final step. 
Similar arguments have already been used in other settings in the literature, often to prove bounds on the location of the maximum of the solution \cite{bogi, lierl, rachh} and also in the context of Hermite-Hadamard inequalities \cite{jianfeng}.

\subsection{Parts of $A$ are close to $x_0$.} The next step is to argue that
$$ \int_{A}   p_t(x_0, y) dy \geq 1 - 4 \varepsilon$$
implies that $A$ has a large intersection with a ball $B_{\sqrt{t}}(x_0)$. We know from the previous section that for $\varepsilon$ sufficiently small the entire ball is contained in $\Omega$. Domain monotonicity of the heat kernel allows to compare the heat kernel $p_t(x_0, \cdot)$ to the heat kernel in $\mathbb{R}^n$ which is strictly larger. This implies
$$ 1 - 4\varepsilon \leq \int_{A}   p_t(x_0, y) dy  \leq \int_{A}  \frac{1}{(4\pi t)^{n/2}} \exp\left( - \frac{\|x_0 -y\|^2}{4t} \right)  dy \leq 1.$$ 
Hence, since the Euclidean heat kernel has total integral 1,
$$  \int_{A^c}  \frac{1}{(4\pi t)^{n/2}} \exp\left( - \frac{\|x_0 -y\|^2}{4t} \right)  dy \leq 4 \varepsilon.$$
Let us now consider $B_{\sqrt{t}}(x_0) \cap A^c$. We clearly have
$$ \int_{ B_{\sqrt{t}}(x_0) \cap A^c}  \frac{1}{(4\pi t)^{n/2}} \exp\left( - \frac{\|x_0 -y\|^2}{4t} \right) \geq  \frac{1}{(4\pi t)^{n/2}} e^{-1/4} \cdot |B_{\sqrt{t}}(x_0) \cap  A^c|.$$
From this we deduce
$$ |B_{\sqrt{t}}(x_0) \cap  A^c| \leq 4 e^{1/4} \varepsilon \cdot (4\pi t)^{n/2} \leq 8\varepsilon  (4\pi t)^{n/2}.$$
We recall that $t = \varepsilon \lambda_k^{-1/2}$ and that
$$ B_{\sqrt{t}}(x_0) = \omega_n t^{n/2},$$
where $\omega_n$ is the volume of the unit ball in $n$ dimensions. Therefore
\begin{align*}
 |B_{\sqrt{t}}(x_0) \cap  A| &=  |B_{\sqrt{t}}(x_0)| -  |B_{\sqrt{t}}(x_0) \cap  A^c| \\
 &\geq \left(1 - \frac{8 \cdot (4\pi)^{n/2}}{\omega_n} \varepsilon\right) |B_{\sqrt{t}}(x_0)|.
 \end{align*}
 Hence there exists a $\varepsilon$ sufficiently small (depending only on the dimension) such that $99\%$ of the volume of $B_{\sqrt{t}}(x_0)$
 is contained in $A$.

 \subsection{Piercing Rays.} We can now conclude the argument as follows (see Fig. \ref{fig:sketch}). The previous section implies that
$$  |B_{\sqrt{t}}(x_0) \cap  A| \geq  \left(1 - \frac{8 \cdot (4\pi)^{n/2}}{\omega_n} \varepsilon\right)|B_{\sqrt{t}}(x_0)|.$$

 \begin{center}
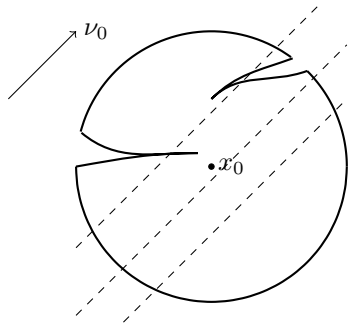
\begin{figure}[h!]
\begin{tikzpicture}[scale=0.9]
\filldraw (0,0) circle (0.04cm);
\node at (0.3, 0) {$x_0$};
\draw [thick,domain=0:45, smooth] plot ({2*cos(\x)}, {2*sin(\x)});
\draw [thick,domain=54:165, smooth] plot ({2*cos(\x)}, {2*sin(\x)});
\draw [thick,domain=180:360, smooth] plot ({2*cos(\x)}, {2*sin(\x)});
\draw[thick] (-2,0) to[out=10, in =180] (-0.2, 0.2) to[out=180, in=320] (-1.92, 0.5);
\draw[thick] (1.4142, 1.4142) to[out=200, in=45] (0,1) to[out=45, in =200] (1.2, 1.6);
\draw[dashed]  (-2+0.7, -2-1+0.7) -- (2,2-1);
\draw[dashed]  (-2, -2+0.8) -- (2-0.4,2+0.4);
\draw[dashed]  (-2, -2-0.2) -- (2,2-0.2);
\draw [->] (-3, 1) -- (-2,2);
\node at (-1.7, 2) {$\nu_0$};
\end{tikzpicture}
\caption{A ball with a small percentage of its mass, the set $A^c$, removed (in practice $A^c$ may look a lot more complicated). There exists a long line that intersects the ball and $A$ along a large set.}
\label{fig:sketch}
\end{figure}
\end{center}

Thus, for $\varepsilon$ sufficiently small, we have that 99\% of the volume of the ball $B_{\sqrt{t}}(x_0)$ is actually contained in $A$.
 We now use the pigeonhole principle (or Fubini's Theorem) to argue that there exists a line in $\mathbb{R}^n$ that points in direction $\nu_0$
 $$ \ell = \left\{ z + t \cdot \nu_0: z \in B_{\sqrt{t}}(x_0) \wedge t \in \mathbb{R} \right\}$$
 such that the line intersects the ball $B_{\sqrt{t}}(x_0)$ along a long segment
$$ | \ell \cap B_{\sqrt{t}}(x_0) | \geq c_4 \sqrt{t},$$
where $c_4 > 0$ depends only on the dimension
 and furthermore
 $$  \frac{| \ell \cap B_{\sqrt{t}}(x_0) \cap A|}{ | \ell \cap B_{\sqrt{t}}(x_0) | } \geq \frac{97}{100}.$$
 \begin{center}
\begin{figure}[h!]
\begin{tikzpicture}[scale=0.7]
\draw [thick,domain=0:360, smooth] plot ({2*cos(\x)}, {2*sin(\x)});
\draw[thick]  (-2.5, 0) -- (0, 2.5);
\draw[thick]  (0, -2.5) -- (2.5, 0);
\draw [->] (-3, 1) -- (-2,2);
\node at (-1.7, 2) {$\nu_0$};
\end{tikzpicture}
\caption{Cutting off the `edges'.}
\label{fig:wings}
\end{figure}
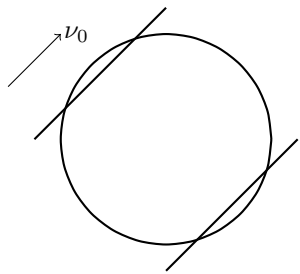
\end{center}
This can be achieved as follows: we choose a small $c_4$ and remove all those parts of a sphere that have intersection less than $c_4 \sqrt{t}$ in direction $\nu_0$ (see Fig. \ref{fig:wings}).
By making $c_4$ sufficiently small (depending only on the dimension), we can ensure that no more than 1\% of the total volume of
the sphere has been lost. We know that at least 99\% of the measure of the ball is in $A$: in the worst case, the `edges' that have
been removed are in $A$ which ensures that at least a $0.98/0.99  \geq 0.98$ portion of the remaining mass is in $A$.
If there were no line going through the remainder of the ball for which
 $$  \frac{| \ell \cap B_{\sqrt{t}}(x_0) \cap A|}{ | \ell \cap B_{\sqrt{t}}(x_0) | } \geq \frac{97}{100},$$
then we can integrate over all the fibers and conclude that at most 97\% of the remaining mass is contained in $A$
which contradicts our inequality.
 
 \subsection{Fundamental Theorem of Calculus.} We conclude with an application of the Fundamental Theorem. Let $\ell$ be the line constructed in the preceding
 section and let $\ell_0 = \ell \cap B_{\sqrt{t}}(x_0)$ denote the line segment that lies fully in the ball $B_{\sqrt{t}}(x_0)$. By construction $\ell_0$ has length at least $| \ell_0| \geq c_4 \sqrt{t}$.
 Denoting the beginning and the endpoint of $\ell_0$ by $a,b$, respectively, we have
 $$ \frac{\partial \phi_k}{\partial \nu_0}(b) -  \frac{\partial \phi_k}{\partial \nu_0}(a) = \int_{\ell_0} \frac{\partial^2 \phi_k}{\partial \nu_0^2} (x) dx.$$
 However, 97\% of the line segment is actually in $A$ and, recalling the definition of
 $$ A =\left\{x \in \Omega:  \frac{\partial^2 }{\partial \nu_0^2} \phi_k(y) \geq  \frac{1}{2}  \frac{\partial^2 \phi_k}{\partial \nu_0^2}(x_0)  \right\},$$
 we have
 $$  \int_{\ell_0} \frac{\partial^2 \phi_k}{\partial \nu_0^2} (x) dx \geq \frac{0.97 \cdot |\ell_0|}{2}  \frac{\partial^2 \phi_k}{\partial \nu_0^2}(x_0)+ \int_{\ell_0 \cap A^c} \frac{\partial^2 \phi_k}{\partial \nu_0^2} (x) dx.$$
 However, we can bound the remaining integral by the supremum norm, meaning 
 $$   \frac{\partial^2 \phi_k}{\partial \nu_0^2}(x) \geq -   \frac{\partial^2 \phi_k}{\partial \nu_0^2}(x_0)$$
  and get that
 $$  \int_{\ell_0 \cap A^c} \frac{\partial^2 \phi_k}{\partial \nu_0^2} (x) dx  \geq - 0.03 \cdot |\ell_0|  \frac{\partial^2 \phi_k}{\partial \nu_0^2}(x_0).$$
 Therefore, recalling that $|\ell_0| \geq c_4 \sqrt{t}$,
 \begin{align*}
 \frac{\partial \phi_k}{\partial \nu}(b) -  \frac{\partial \phi_k}{\partial \nu}(a)  &\geq 0.45 \cdot |\ell_0|   \frac{\partial^2 \phi_k}{\partial \nu_0^2}(x_0) \\
 &\geq 0.45 \cdot c_4 \cdot \sqrt{\varepsilon} \cdot \lambda_k^{-1/2}   \frac{\partial^2 \phi_k}{\partial \nu_0^2}(x_0).
 \end{align*}
 We combine this with the  upper bound
 $$ \frac{\partial \phi_k}{\partial \nu}(b) -  \frac{\partial \phi_k}{\partial \nu}(a)   \leq 2 \cdot \| \nabla \phi_k\|_{L^{\infty}} \leq 2 c_2 \cdot \lambda_k^{\frac{n+1}{4}}$$
 and, since $\varepsilon$ is an absolute constant, this is the desired result.

 \subsection{Concluding Remarks.} It seems that certain variations of this argument are feasible and some may prove to be useful in other settings. We quickly highlight one such variation. In the first half of the proof, we establish that if there is a large second derivative, then there are large second derivatives in a large subset of the $\sim \sqrt{t}$ neighborhood around the point. If the Hessian is very large in a point, then it is also large in a large set in the neighborhood. Bochner's formula states that
 $$ \frac{1}{2} \Delta |\nabla u|^2 = \left\langle \nabla \Delta u, \nabla u \right\rangle + \|D^2 u \|^2.$$
 In the case of a Laplacian eigenfunction, the identity simplifies to
  $$ \frac{1}{2} \Delta |\nabla \phi_k|^2 = - \lambda_k \| \nabla \phi_k\|^2 + \|D^2 \phi_k \|^2.$$
If the maximum size of the Hessian exceeds $\sqrt{\lambda_k} \|\Delta \phi_k\|_{L^{\infty}}$ by a large constant in a point $x_0 \in \Omega$, then
  $$  \Delta |\nabla \phi_k|^2 \sim \|D^2 \phi_k \|^2 \sim \lambda_k^{\frac{n+3}{2}} \qquad \qquad (\diamond)$$
This, however, shows that $|\nabla \phi_k|^2$ is undergoing growth. We remark that if an arbitrary $f:\mathbb{R}^n \rightarrow \mathbb{R}$ is nonnegative $f \geq 0$
and satisfies an inequality of the type
$$ \Delta f \geq C \qquad \mbox{in}~B_R(0), \qquad \mbox{then} \qquad \|f\|_{L^{\infty}(B_R(0))} \gtrsim_n C \cdot R^2.$$
This can be seen by comparing $f$ with the function $g$
$$ g(x) = \frac{C}{2} \frac{\|x\|^2}{2n}.$$
We have $\Delta (f-g) \geq C/2$ and thus $f-g$ assumes its maximum on the boundary. We also know that $f(0) - g(0) = 0$ and thus the maximum is nonnegative and thus at least of size $\sim CR^2$. Applying this to our inequality $(\diamond)$
in a $B_{\sqrt{\varepsilon} \lambda_k^{-1/2}}(x_0)$ neighborhood, we obtain
$$ \left\|  \nabla \phi_k \right\|^2_{L^{\infty}\left(B_{\sqrt{\varepsilon} \lambda_k^{-1/2}}(x_0)\right)} \gtrsim_{\varepsilon, n} \lambda_k^{\frac{n+3}{2}} \lambda_k^{-1} = \lambda_k^{\frac{n+1}{2}} .$$
For sufficiently large $C$, this then contradicts the gradient estimate. A difficulty with this argument is that $(\diamond)$ is not a priori satisfied in the entire ball but merely in a very large subset: however, in the complement of that set, we always have
  $$ \frac{1}{2} \Delta |\nabla \phi_k|^2 = - \lambda_k \| \nabla \phi_k\|^2 + \|D^2 \phi_k \|^2 \geq - \lambda_k \| \nabla \phi_k\|^2$$
which cannot be arbitrarily negative and is a constant factor smaller than the Hessian. One would expect that a more refined maximum principle could then be applied in that case. This type of argument may be simpler to apply on a manifold since it does not rely on arguments along sub-manifolds.

\end{document}